%Paper: alg-geom/9308001
%From: ALBANO@dm.unito.it
%Date: Tue,  3 AUG 93 16:45 GMT

% This is a AmSTeX 2.1 file

\documentstyle{amsppt}
\refstyle{C}

\mag1200

\pagewidth{13cm}
\pageheight{19cm}

%%%% general definitions  %%%%%%%

%% spacing
\def\oneline{\vskip12pt}

%% symbols of frequent use
\def\cy{\lambda}                  % \cy  = primitive cycle on X
\def\Hy{\Sigma}                   % \Hy  = hyperplane section of X

\def\Z{\Bbb Z}                    % \Z   = relative numbers
\def\Q{\Bbb Q}                    % \Q   = rational numbers
\def\C{\Bbb C}                    % \C   = complex numbers
\def\O{\Cal O}                    % \O   = structure sheaf
\def\P{\Bbb P}                    % \P   = projective space

      % \Aut = automorphisms group
\def\Gal{\operatorname{Gal}}      % \Gal = Galois group
\def\Res{\operatorname{Res}}      % \Res = residue operator
\def\GM{\overline{\nabla}}        % \GM  = Gauss-Manin connection
\def\mi#1{\bold\mu_{#1}}              % \mi  = roots of unity
\def\cp{\bullet}                  % \cp  = cup-product
\def\GQ#1{G^d_{\Q}(#1)}           % \GQ  = Griffths group tensored by Q

 %%%% end general definitions %%%%

\topmatter

\title
On the Griffiths group of the cubic sevenfold
\endtitle

\author  A. ALBANO -- A. COLLINO  \endauthor

\affil Universit\`a di Torino \endaffil

\address Alberto Albano, Dipartimento di Matematica, Universit\`a di Torino,
Via Carlo Alberto~10, 10123~Torino, ITALY  \endaddress

\email  albano\@dm.unito.it
\endemail

\address Alberto Collino, Dipartimento di Matematica, Universit\`a di Torino,
Via Carlo Alberto~10, 10123~Torino, ITALY  \endaddress

\email  collino\@dm.unito.it
\endemail

\date
2 August 1993
\enddate

\thanks Partially supported by Science Project ``Geometry of Algebraic
Varieties'', n.~0-198-SC1, and by fundings from M.U.R.S.T. and
G.N.S.A.G.A.~(C.N.R.), Italy  \endthanks

\subjclass 14C30, 14C10, 14J40, 14K30  \endsubjclass

\abstract  In this paper we prove that the Griffiths group of a general
cubic sevenfold is not finitely generated, even when tensored with~$\Q$.
Using this result and a theorem of Nori, we
provide examples of varieties which have some Griffiths group not
finitely generated but whose corresponding intermediate Jacobian is
trivial.
\endabstract

\endtopmatter

\document

\head 1. Introduction  \endhead

For $X$ a smooth complex projective variety, one of the important questions in
algebraic geometry is the existence and the behavior of non trivial
subvariety on~$X$, i.e., subvarieties that are not the intersection of~$X$
with a hypersurface. In general, the approach is as follows: we let $Z^d(X)$ be
the free abelian group generated by the irreducible subvarieties of
codimension~$d$ in~$X$, and introduce various equivalence relations
on~$Z^d(X)$.
The problem is then to understand the quotients of $Z^d(X)$ under these
relations. The first, and most obvious, is homological equivalence, but there
are other interesting relations, namely algebraic and rational equivalence. The
so--called classical cases, when $d=1$ (divisors) or $d=\dim X$ ($0$--cycles)
are well understood, but little is known in general. In these cases,
homological and algebraic equivalence coincide, and the group of cycles
algebraically equivalent to zero modulo rational equivalence is isomorphic
to an
abelian variety.

One of the first result for general~$d$ is
due to Griffiths \cite{8}, who showed that  homological and algebraic
equivalence do not coincide if $d\ne1, \dim X$.
The group of cycles on~$X$ of codimension~$d$ homologous to zero modulo
algebraic  equivalence is called the $d$-th~Griffiths group of~$X$ and is
denoted by~$G^d(X)$; we will mainly consider the $\Q$--vector space
$\GQ X = G^d(X)\otimes\Q$, or equivalently, the Griffiths group modulo
torsion.
Clemens \cite{5} then showed
that $\GQ X$ can have infinite dimension over~$\Q$.
Clemens' example is given by a general quintic hypersurface in~$\P^4$, and this
opened the way to the construction of other examples (see \cite{1}, \cite{2},
\cite{13}, \cite{17}).

In each of these examples, the techniques of proof were quite ad~hoc, and
there seemed not to be a general pattern. Only recently, two fundamental
results
have been obtained: the first is Voisin's proof of Clemens' theorem on the
quintic threefolds \cite{21}, and the other is Nori's beautiful theorems on
Hodge theoretic connectivity \cite{14}.

In showing that $\GQ X$, for some variety~$X$, is infinite dimensional,
there are basically two steps: first, one must find infinitely many cycles and
then one has to show that they are independent. In all of the previous proofs,
the cycles were found by geometrical considerations, which cannot generalize to
more complicated situations. Voisin gives a very general method for finding
cycles, based on an argument of Green concerning the Noether--Lefschetz locus
of a suitably defined family of varieties; then, using Griffiths' description
of
cohomology via residues of rational differential forms, she is able to compute
explicitely (this is a highly non trivial computation) an infinitesimal
invariant
attached to an algebraic cycle, which permits to conclude that the cycles so
found are independent.
In~\cite{21}, she
works out the case of quintic threefolds, but her technique is quite general
(Bardelli and M\"uller-Stach have worked out the case of the complete
intersection of two cubics in~$\P^5$).

A result in the opposite direction is Green's theorem  (\cite{7} and Voisin,
unpublished):  other than for cubic hypersurfaces of dimension~3 and~5 and
for quartic hypersurfaces of dimension~3, for which the intermediate Jacobian
is
an abelian variety and the Abel--Jacobi map is known to be surjective,
the only  cases in which the image of
the Abel--Jacobi map for a general hypersurface may be non~torsion , i.e.,
$\GQ X$ may be non zero,
are quintic
threefolds (where the image is  non~torsion [5]) and  cubic sevenfolds,
i.e., cubic hypersurfaces in~$\P^8$.   Moreover, the middle cohomology of a
smooth cubic sevenfold carries a Hodge structure of level~3, similar to the one
of a smooth quintic threefold.
As Voisin's proof is Hodge--theoretic, we were
led to see if her techniques could be carried over to this case: in fact,
following \cite{21} we prove the following:

\proclaim{Theorem 1} For a general smooth cubic sevenfold~$X$, the group of
$3$--cycles homologous to zero modulo algebraic equivalence is not finitely
generated, even when tensored with $\Q$.
\endproclaim

The main difference with the previously known cases of infinite
generation is that the cycles are of codimension~four and they are not
curves. This makes it more difficult to treat the situation geometrically, and
in fact we use in essential way Voisin's method of doing things
Hodge--theoretically to overcome this difficulty.

On the other hand, our example is not too surprising:
Voisin has recently shown  (\cite{22}) that the Griffiths
group is non~torsion for any Calabi--Yau threefold, i.e., a K\"ahler
threefold~$X$
such that $K_X$ is trivial and $h^{1,0}(X) = h^{2,0}(X) = 0$.
It is an open question whether the Griffiths group is always infinitely
generated for such threefolds.
Calabi--Yau threefolds have attracted much interest from physicist, in
connection
with problems in string theory. One of the main conjectures (based on physics
considerations) is the so--called ``mirror symmetry'', that roughly states that
to each Calaby--Yau threefold~$X$ (plus some extra data) is attached another
Calabi--Yau threefold~$\tilde X$, the ``mirror of~$X$'', whose Hodge diamond
is the mirror of the Hodge diamond of~$X$. This conjecture is full of
mathematical implications, some of which have been verified (see \cite{12} for
further information, and for a precise statement of the mirror conjecture).
One potential counterexample to
this conjecture is given by rigid threefolds, i.e.,  with $h^{2,1}(X)=0$, since
then the mirror should have $h^{1,1} = 0$ and so cannot be K\"ahler. In a
recent preprint \cite{3}, Candelas, Derrick and Parkes
 show how a cubic sevenfold can be
considered the mirror (in an appropriate sense) of one such rigid Calabi--Yau
threefold. So, in some sense, the cubic sevenfold should have properties
similar to those of Calabi--Yau threefolds: Theorem~1 shows that this is so
with respect to the Griffiths group, which is not only non torsion but even
infinitely generated (modulo torsion).

But things can be even worse.  The essential
tool that has always been used to show that two cycles are algebraically
independent is the Abel--Jacobi map~$\theta_d$ from the group of cycles of
codimension~$d$ homologous to zero modulo rational equivalence to the Griffiths
intermediate Jacobian~$J^d(X)$.  In~\cite{14}, Nori introduces a deeper
version of Griffiths approach: his theorems allow to show independence even
when
the Abel--Jacobi map is identically zero. He provides the first
example of cycles homologically equivalent to zero but not algebraically
equivalent to zero on a particular variety without the use of the Abel--Jacobi
map on the variety, and in fact belonging to the kernel of the Abel--Jacobi
map. He shows that the algebraic part of the primitive
cohomology of a hypersurface~$X$ maps injectively under restriction to a
subgroup of the Griffiths group of a complete intersection of~$X$ with
general hypersurfaces of sufficiently high degree. However, in this way
one cannot obtain infinitely generated subgroups of the Griffiths group.

To get such examples, we use Theorem~1 and Nori's theorem. Let $X$ be
a general smooth cubic sevenfold, $D_1$,
$D_2$ two general hypersurfaces of sufficiently high degree and let $Y=X\cap
D_1\cap D_2$.

\proclaim{Theorem 2} The (rational) Griffiths group $G^4_{\Q}(Y)$ of
codimension~$4$ cycles
in~$Y$ is not finitely generated even if the intermediate Jacobian~$J^4(Y)$
is zero.
\endproclaim

The same proof shows
that also for $Y' = X\cap D_1$ the Griffiths group
is not finitely generated, so we have examples in both even and
odd dimension.

We conclude this introduction explaining in details how Theorem~2 follows from
Nori's theorem and Theorem~1.
Let $X$ be a
smooth projective variety~$X$,
and let $CH^d(X)$ be the Chow group of cycles of codimension~$d$
on~$X$: Nori defines an increasing filtration
$A_{\bullet}CH^d(X)$, whose main property is that if $\eta\in
A_rCH^d(X)\otimes\Q$ then there is an algebraic subset~$Z\subset X$ of pure
codimension~$(d-r-1)$ so that $Z$ contains the support of~$\eta$, and the
fundamental class of~$\eta$ vanishes in $H_{2n-2d}(Z,\Q)$, where $n=\dim X$. In
particular, $A_0CH^d(X)$ consists of cycles algebraically equivalent to zero.

If $\eta\in A_rCH^d(X)$, then not only $\eta$ is homologous to zero, but
$\theta_d(\eta)$ has to satisfy some restriction:
if $J_r^{d}(X)$ is the intermediate
Jacobian of the largest integral Hodge structure contained in
$F^{d-r-1}H^{2d-1}(X)$, then $\theta_d(\eta)\in J_r^{d}(X)$. Note that
$J^{d}_0(X)$ is the ``largest
abelian subvariety'' of the intermediate Jacobian~$J^{d}(X)$
We can now
state \cite{14, Theorem~2}:

\proclaim{Theorem} {\rm(Nori)}\ Let $X \subset \P^m$ be smooth, projective,
and let $Y$ be the intersection of~$X$ with $h$~general hypersurfaces of
sufficiently large degrees. Let $\xi\in CH^d(X)$ and let $\eta = \xi|Y$.
Assume that $r+d<$~dimension of~$Y$. If $\eta\in A_rCH^d(Y)\otimes\Q$, then
\roster
\item The cohomology class of $\xi$ vanishes in $H^{2d}(X, \Q)$, and
\item the Abel--Jacobi image of a non--zero multiple of $\xi$ belongs to
$J_r^{d}(X)$.
\endroster
\endproclaim

It is now clear how to prove Theorem~2:  Let
$R: CH^4_{\text{hom}}(X) \to CH^4_{\text{hom}}(Y)$ be the restriction map
and $\pi: CH^4_{\text{hom}}(Y) \to G^4(Y)$ be the canonical projection.
Let $\xi\in CH^4_{\text{hom}}(X)$:  if $\pi\circ R(\xi)=0$ then, by
Nori's theorem (with $r=0$), we would have
$\theta_4(\xi)\in J_0^{4}(X)$. But for
general~$X$, a monodromy argument shows that $J_0^{4}(X)= 0$ and hence
$\ker(\pi\circ R) \subseteq \ker\theta_4$.
Since
for a general cubic sevenfold~$X$ we
have, from the proof of Theorem~1, that the image
$\theta_4(CH^4_{\text{hom}}(X))$ is not finitely generated,
even when tensored with~$\Q$, the same is
true for $G^4(Y)\otimes\Q$.

\head 2. Cycles on the cubic sevenfold \endhead

In this and the following sections we show how Voisin's argument can be used
to prove that
the Griffiths group of a generic cubic sevenfold is infinitely generated.
There are three parts in the proof: the first two (\S2 and 3)
are very similar to the
corresponding ones in \cite{21}, and so  we indicate only
the changes needed in the present case, giving precise references for all
the results we quote. The computation of the infinitesimal invariants in the
last part  however (\S4), is
quite different, and we will give full details. For ease of reference, we use
the same notations as in \cite{21}.

Let $(\cy, \Hy, X)$ be a triple, where
$X \subset \P^8$ is a smooth cubic hypersurface
defined by an equation~$F$, $\Hy \subset X$
is a smooth hyperplane section, and $\cy \in H^6(\Hy, \Z)^{\text{prim}} \cap
H^{3,3}(\Hy)$ is the class of a primitive algebraic 3--cycle on~$\Hy$.

We note the following two useful facts.

\proclaim{Lemma 2.1} Let $X$ and $\Hy$ be as above. Then:

\roster
\item $h^{7,0}(X) = h^{6,1}(X) = 0$; $h^{5,2}(X) = 1$; $h^{4,3}(X) = 84$.
\item $h^{6,0}(\Hy) = h^{5,1}(\Hy) = 0$; $h^{4,2}(\Hy) = 8$;
      $h^{3,3}(\Hy) = 36$;
\item $h^0(N_{\Hy/X}) = 8$.
\endroster
\endproclaim

\demo{Proof} The simplest way to compute these Hodge numbers is to
use Griffiths' description of the Hodge filtration of~$X$
and~$\Hy$ with rational forms with poles along the hypersurface \cite{8}.
Since the Hodge numbers are invariant under deformations,
one can compute them for the Fermat hypersurface, whose Jacobian ideal is
particularly simple.
Alternatively, one can use the formulas for the Hodge numbers of complete
intersections given in \cite{6, Th\'eor\`eme 2.3}.

The last formula is clear.
\qed
\enddemo

Hence, the Hodge numbers of~$X$ and~$\Hy$ ``look like'' the ones of a
quintic threefold and its hyperplane section and moreover $h^{4,2}(\Hy) =
h^0(\Hy, N_{\Hy/X})$.

The second observation is crucial.

\proclaim{Proposition 2.2} The Hodge $(3,3)$--conjecture is true for any
smooth cubic sixfold $\Hy \subset \P^7$.
\endproclaim

\demo{Proof} This is due to Steenbrink    \cite{19}.\qed
\enddemo

\oneline

Let $\Cal G$ be the smooth variety that parametrizes the pairs~$(X,\Hy)$
such that $\Hy \subset X$, and let $\cy \in H^6(\Hy, \Z)^{\text{prim}} \cap
H^{3,3}(\Hy)$. Since $\cy$ is real,
by Lemma 2.1, the condition ``$\cy$ remains of type~(3,3) in
$H^6(\Hy)$'' is given by $h^{4,2}(\Hy) = \dim H^2(\Hy,
\Omega_{\Hy}^4)$ equations on $\Cal G$, so that locally the
family~$\Cal F_{\cy}$ of the triples $(\cy, \Hy, X)$ is a
subvariety of~$\Cal G$ of codimension at most $h^{4,2}$. In
particular, $\Cal F_{\cy}$ is smooth at any point in which its
tangent space is of codimension~$h^{4,2}$ in the tangent space
of~$\Cal G$.

The exact sequence
$$
\CD
0 @>>> T_{\Hy} @>>> T_{X|\Hy}  @>>> \O_{\Hy}(1)  @>>> 0
\endCD   \tag1
$$
induces a Kodaira--Spencer map $\rho : H^0(\O_{\Hy}(1)) \to H^1(T_{\Hy})$.
We assume that $\rho$ is injective (this will be true in our situation).

Let $\alpha : T_X \to T_{X|\Hy}$ be the restriction map, and
$\beta : T_{\Hy} \to T_{X|\Hy}$  the inclusion map.
Also, we denote by~$\cp$  the cup--product
$H^1(\Hy, T_{\Hy}) \otimes H^3(\Omega_{\Hy}^3) \to
H^4(\Omega_{\Hy}^2)$.

With these notations, we have:

\proclaim{Lemma 2.3} Assume that $\rho$ is injective. Then the infinitesimal
deformations of the triple $(\cy, \Hy, X)$ are parametrized by the subspace
$T{\Cal F}_{\cy(\Hy, X)} \subset H^1(T_{\Hy}) \times H^1(T_X)$ given by:
$$
(u, v) \in T{\Cal F}_{\cy(\Hy, X)} \iff \alpha(v) = \beta(u) \text{ \rm in }
H^1(T_{X|\Hy}) \text{ \rm and } u\cp\cy^{3,3} = 0 \text{ \rm in }
H^4(\Omega^2_{\Hy}) $$
\endproclaim

\demo{Proof} \cite{21, 1.1 and 1.3} \qed
\enddemo

Let $\pi : T{\Cal F}_{\cy(\Hy, X)} \to H^1(T_X)$ be the map induced by the
second projection. Consider the diagram induced by~(1):
$$
\CD
H^0(\O_{\Hy}(1))  @>{\rho}>>  H^1(T_{\Hy})  @>{\beta}>> H^1(T_{X|\Hy}) @>>>
 H^1(\O_{\Hy}(1))  \\
@.                                @VV{\cy^{3,3}}V    \\
                  @.          H^4(\Omega^2_{\Hy})    \\
\endCD  \tag2
$$

Since $H^1(\O_{\Hy}(1)) = 0$, a simple diagram chase shows that $\pi$ is
surjective if and only if the maps $\cy^{3,3}$ and $\cy^{3,3} \circ \rho$
have the same image. Moreover, if $\rho$ is injective, then $\pi$ is
injective if and only if $\cy^{3,3}\circ\rho$ is injective.

Now, $\dim H^0(\O_{\Hy}(1)) = \dim H^4(\Omega^2_{\Hy})$, and so
$\cy^{3,3} \circ \rho$ is injective if and only if it is surjective, and
in this case $\pi$ is an isomorphism, i.e., the family $\Cal F_{\cy}$
project isomorphically onto the family~$\Cal Y$ of deformations of~$X$.

The description of the map $\cy^{3,3}\circ\rho$ in
terms of polynomial multiplication in the Jacobian ring of~$\Hy$ is given
in \cite{21, 1.4}. Since of course the grading is different in our case, we
repeat it here, mainly to fix notations.

Let $(x_0, \dots, x_8)$ be homogeneous coordinates in~$\P^8$ such that
$\Hy = X \cap \{x_8 = 0\}$. Let $F$ be the equation of~$X$ and let
$G = F_{|\P^7}$ be the equation of~$\Hy$. The Jacobian ring~$R(\Hy)$ is defined
as:
$$
R(\Hy) = \C[x_0, \dots, x_7]/ \left( \frac{\partial G}{\partial x_0}, \dots,
  \frac{\partial G}{\partial x_7} \right)
       = S(\P^7) / J(G).
$$
$R(\Hy)$ is a graded ring, and by Griffiths' residue theory there are
well--known isomorphisms
\roster
\item $H^1(T_{\Hy}) \cong R^3(\Hy)$,
\item $H^3(\Omega^3)^{\text{prim}} \cong R^4(\Hy)$,
\item $H^4(\Omega^2) \cong R^7(\Hy)$,
\item $H^1(T_{\Hy}(-1)) \cong R^2(\Hy)$.
\endroster

These isomorphisms are compatible in the following sense: if a polynomial
class $P_{\cy}\nobreak\in\nobreak R^4(\Hy)$ corresponds under~(2) to
a primitive cohomology class~$\cy^{3,3}$, then the map
$\cy^{3,3}:H^1(T_{\Hy})\to H^4(\Omega^2)$ corresponds, under the
isomorphisms~(1) and~(3), to multiplication  by~$P_{\cy}$.

Let now $e \in R^2(\Hy)$ be the image of the
polynomial~$(\partial F/\partial x_8)_{|\P^7}$. Under the isomorphism~(4),
$e$~corresponds to the extension class of the exact sequence
$$
\CD
0 @>>> T_{\Hy} @>>> T_{X|\Hy}  @>>> \O_{\Hy}(1)  @>>> 0
\endCD
$$
Since the isomorphisms (1) and (4) are compatible with polynomial
multiplication, we obtain that $\rho : H^0(\O_{\Hy}(1)) \to H^1(T_{\Hy})$
is given by multiplication by~$e$ from $R^1(\Hy) \cong H^0(\O_{\Hy}(1))$
to $R^3(\Hy)$.

{}From the above discussion we obtain:

\proclaim{Proposition 2.4} Let $(\cy, \Hy, X)$ be a triple as above. If
the multiplication by $P_{\cy}\nobreak\cdot\nobreak e$ induces an isomorphism
$R^1(\Hy) \cong R^7(\Hy)$, then the map $\pi : T{\Cal F}_{\cy} \to H^1(T_X)$
is an isomorphism, and hence the family~$\Cal F_{\cy}$ of deformations
of the triple~$(\cy, \Hy, X)$ is smooth, of codimension~$h^{4,2}(\Hy)$
in~$\Cal G$, and projects isomorphically (locally) onto the family~$\Cal Y$
of deformations of~$X$.
\endproclaim

This is \cite{21, Proposition 1.5}. We only note that the hypothesis of the
proposition implies that the Kodaira--Spencer map is injective, and then we can
use Lemma~2.3.

\medskip
Let ${\Cal S}(X)$ be the Noether--Lefschetz locus for the hyperplane sections
of~$X$, i.e., the set of those~$\Hy$ such that $H^6(\Hy, \Z)^{\text{prim}} \cap
H^{3,3}(\Hy) \ne 0$. If $\pi$ is an isomorphism at a point $(\cy, \Hy, X)$,
then
the component ${\Cal S}_{\cy}$ of ${\Cal S}(X)$ given by $\cy$, in a
neighborhood of $\Hy$ consists only of the isolated reduced point~$\Hy$, and
has the expected codimension equal to $h^{4,2}(\Hy)$. Hence, as in \cite{21,
1.6}, we can use Green's argument (\cite{4, \S5}) to obtain:

\proclaim{Proposition 2.5} If there exists a triple $(\cy, \Hy, X)$ such that
$\pi : T{\Cal F}_{\cy} \to H^1(T_X)$ is an isomorphism, then the
$0$--dimensional components of the Noether--Lefschetz locus ${\Cal S}(X)$ form
a countable dense subset of $\P(H^0(\O_X(1)))$.
\endproclaim

We will give later (Proposition~4.2)
an explicit example of such a triple, using Proposition~2.4.
This will show that there is a Zariski open set of the family~$\Cal Y$
of deformations of~$X$ over which such
triples exist.

\head 3. The infinitesimal invariant of normal functions  \endhead

To each triple $(\cy, \Hy, X)$ as above, we can associate an element
$Z_{\cy} \in CH_3^{\text{Hom}}(X)$, the group of 3--cycles on~$X$
homologous to~0, modulo rational equivalence, defining  $Z_{\cy} = j_*\cy$,
where $j : \Hy \to X$ is the inclusion map. In fact, by Proposition~2.2, $\cy$
gives a primitive algebraic 3--cycle on~$\Hy$, and hence $Z_{\cy}$ is
homologous to~0 in~$X$.

If $(\cy, \Hy, X)$ satisfies the hypothesis of Proposition~2.4, the local
family~$\Cal F_{\cy}$ is isomorphic to~$\Cal Y$ in a neighborhood of the
point $o \in \Cal Y$ corresponding to~$X$. Hence there is a universal
family $p : \Cal X \to \Cal Y$ and a relative cycle $\Cal Z_{\cy} \in
CH_3^{\text{Hom}}(\Cal X/\Cal Y)$. From this one obtains a normal
function~$\nu_{\cy}$, with infinitesimal invariant~$\delta\nu_{\cy}$ (we will
recall later the definition of~$\delta\nu_{\cy}$).

Let $q : \Cal J \to \Cal Y$ be the intermediate Jacobian bundle  whose fibre in
$t\in\Cal Y$ is $q^{-1}(t) = J({\Cal X}_t)$, the intermediate Jacobian
of~${\Cal X}_t$.
The normal function $\nu_{\cy}$ is the holomorphic section of this bundle
given by $\nu_{\cy}(t) = \Phi_{\Cal X_t}(\Cal Z_{\cy}(t))$, where
$\Phi_{\Cal X_t}$ is the Abel--Jacobi map of ${\Cal X}_t$.

By definition,
$\Cal Z_{\cy} = \Cal Z^1_{\cy} - \Cal Z^2_{\cy}$, where  for~$t$ in a small
disk around~$0$, $\Cal Z^i_{\cy}(t)$ is contained in a hyperplane section
$\Hy_t \subset \Cal X_t$, varying holomorphically with~$t$.
It would be convenient to have $\Cal Z^i_{\cy}(t)$ smooth in~$\Hy_t$, but
this is not necessarily true. However, by a result of Kleiman \cite{11},
it is possible to find a smooth cycle rationally equivalent to $2(\Cal
Z^i_{\cy}(t) + D(t))$, where $D(t)$ is a sufficiently large multiple of a
hyperplane section of codimension~3 in $\Hy_t$. This simply means that we have
to multiply  by~2 all the formulas in~\cite{21}, and hence we might as well
assume that the cycles $\Cal Z^i_{\cy}(t)$ are already smooth.
We can also, up to
rational equivalence inside $\Hy_t$, assume that the intersection of
$\Cal Z^1_{\cy}(t)$ and $ \Cal Z^2_{\cy}(t)$ is transversal and consists of
finitely many points. Blowing up these points, we can assume that
$\Cal Z^1_{\cy}(t)$ and $\Cal Z^2_{\cy}(t)$ do not meet inside $\Hy_t$ and
are homologous in~$\Cal X_t$. We now choose a differentiable chain
$\Gamma_t$, varying smoothly with~$t$, such that $\partial\Gamma_t =
\Cal Z^1_{\cy}(t) - \Cal Z^2_{\cy}(t)$.

Let $\Cal H^7 = R^7p_*\C\otimes\Cal O_{\Cal Y}$, $\Cal H^{i,j} =
R^jp_*(\Omega^i_{\Cal X/\Cal Y})$ and $F^i\Cal H^7 = \oplus_{k\ge i}\Cal
H^{k, 7-k}$. Then the normal function~$\nu_{\cy}$ lifts to a holomorphic
section~$\psi_{\cy}$ of $(F^4\Cal H^7)^{\vee}    $ by
$$
\psi_{\cy}(t)(\omega_t) = \int_{\Gamma_t} \omega_t
$$
where $(\omega_t)_{t\in\Cal Y}$ is a section of $(F^4\Cal H^7)$.

We briefly recall here  Griffiths' definition of the infinitesimal
invariant~$\delta\nu_{\cy}$: let
$$
\GM : \Cal H^{4,3}\otimes T_{\Cal Y} \to \Cal H^{3,4}
$$
be the variation of Hodge structure on $\Cal Y$. In our case $\GM$
is surjective in every point and so $\ker\overline\nabla$ is a vector bundle
over~$\Cal Y$. The infinitesimal invariant~$\delta\nu_{\cy}$ is the section
of~$(\ker\GM)^{\vee}$ defined in the following way: let
$\sum_i\omega_i\otimes\chi_i \in (\ker\GM)_{t_0}$ and let $\tilde\omega_i(t)$
be a section of $F^4\Cal H^7$ in a neighborhood of~$t_0$, such that
$\tilde\omega_i(t_0) = \omega_i$. By definition of $\ker\GM$,
$\sum_i\nabla_{\chi_i}(\tilde\omega_i) \in F^4\Cal H^7(t_0)$, and hence we can
define $$
\delta\nu_{\cy}(\sum_i\omega_i\otimes\chi_i) = \sum_i
\chi_i(\psi_{\cy}(\tilde\omega_i)) - \psi_{\cy}(t_0)(
\sum_i\nabla_{\chi_i}(\tilde\omega_i)).
$$

By quasi--horizontality of normal functions, this formula gives a well defined
element of $(\ker\GM)_{t_0}^{\vee}$, i.e., it does not depend on the choice of
the~$\tilde\omega_i$ \cite{9}.

We can now repeat the argument of \cite{21, 2.4--2.12} to rewrite the
infinitesimal invariant in terms of Jacobian rings. We explicitely note two
facts: first of all, the proof of
the crucial Proposition~2.4 of \cite{21} holds in our case
(with an obvious change in the grading and dimension) since the cycles $\Cal
Z_{\cy}$ are smooth. Secondly,
if $R(X)$ is the jacobian ring of $X$, let $\mu_X : R^3(X)\otimes R^3(X) \to
R^6(X)$ be the multiplication map. Then, as in \cite{21, 2.5}, we can identify
the map $\GM_{(0)}$ with $\mu_X$.

The formula \cite{21, 2.4} holds only for tensors of rank~one. To obtain
a general formula we need:

\proclaim{Lemma 3.1} For $X$ generic, $\ker \mu_X$ is generated by tensors of
rank one.
\endproclaim

\demo{Proof} A simple proof can be given in the following way:
let $X$ be the
Fermat hypersurface of degree~$3$; then the Jacobian ideal of~$X$ is generated
by $x_0^2, \dots, x_n^2$, so that a basis for $R^3(X)$ is given by the
classes of the monomials
$x_ix_jx_k$ with $i<j<k$.

Since $(x_ix_j)x_k\otimes x_a(x_bx_c) - (x_ix_j)x_a\otimes x_k(x_bx_c)  =
x_ix_j(x_a+x_k)\otimes x_bx_c(x_a-x_k) - x_ix_jx_a\otimes x_ax_bx_c +
x_ix_jx_k\otimes x_bx_cx_k$, we see that
any monomial $P\otimes Q = x_ix_jx_k\otimes x_ax_bx_c$
can be rewritten, modulo tensors of rank~1 belonging to~$\ker\mu_X$, with
the indices $\{i,j,k,a,b,c\}$ in strictly increasing order (or else it
already belonged to $\ker\mu_X$).

Let $\omega =
\sum_iP_i\otimes Q_i \in \ker\mu_X$. The observation above allows
to rewrite, modulo tensors of rank~$1$
belonging to $\ker\mu_X$,
all the summands in a standard form that must then be zero.

To prove the claim for~$X$
generic, one  observes that $\mu_X$ induces a linear map $\P(R^3\otimes
R^3)\to\P(R^6)$ and the tensors of rank one are simply the image of
$\P(R^3)\times\P(R^3)$ in $\P(R^3\otimes R^3)$ via the Segre map. For $X$ the
Fermat hypersurface, the intersection of this image with $\P(\ker\mu_X)$
generates  $\P(\ker\mu_X)$ and so the same will be true generically.
\qed
\enddemo

{}From the this discussion, one
obtains the following formula for
the infinitesimal invariant $\delta\nu_{\cy}$
in terms of Jacobian ideals:
Let $\mu_{\Hy} : R^3(\Hy) \otimes R^3(\Hy) \to R^6(\Hy)$ be the multiplication
map. Then, as Voisin shows in an ingenious way in \cite{21, 2.11 and 2.12},
the infinitesimal invariant $\delta\nu_{\cy}$ can
be viewed as an element of $(\ker \mu_{\Hy})^{\vee}$, i.e., it can be computed
using the hyperplane section~$\Hy$ of~$X$. With the same argument we obtain

\proclaim{Proposition 3.2} If
the multiplication by $P_{\cy} \cdot e$ induces an isomorphism
$f_{\cy} : R^1(\Hy) \to R^7(\Hy)$, then
$$
   \delta\nu_{\cy}(\sum_i P_i\otimes R_i) = \sum_i P_{\cy} P_i (f_{\cy}^{-1}
(P_{\cy}R_i)) \in R^8(\Hy) \cong \C.
$$
\endproclaim
\demo{Proof} \cite{21, 2.12}  \qed
\enddemo

\head 4. Infinite generation of the Griffiths group  \endhead

We can now study the independence over~$\Q$ of the $\delta\nu_{\cy}$. For this,
we apply the formula of proposition~3.2  to the cohomology classes $\cy$
of a cubic sixfold~$\Hy$ for which the rank of $H^3(\Omega^3_\Hy)^{\text{prim}}
\cap H^6(\Hy, \Z)$ is at least~2. One such is
the Fermat cubic sixfold~$\Hy$, whose equation is
$F = x_0^3 + \dots + x_7^3$. The structure of the Hodge classes of~$\Hy$ is
well known, as it is for all Fermat hypersurfaces. We recall here the
relevant facts for the case at hand. Complete proofs, in the general case,
can be found in \cite{15}, \cite{16}, \cite{18}. We also note again that since
the Hodge conjecture is true for $\Hy$, it is enough to work with rational
cohomology classes of type~$(3,3)$, since these will correspond to algebraic
$3$--cycles on~$\Hy$.

Let $\mi3$ be the group of third roots of unity, let $\zeta = e^{2\pi i/3}$
and let
$G = (\mi3)^8/\Delta$, where $\Delta$ is the diagonal. $G$ acts on $\Hy$ by
$$
(\zeta_0, \dots, \zeta_7) \cdot (x_0 : \dots : x_7) = (\zeta_0 x_0: \dots :
\zeta_7 x_7).
$$
The group of characters of $G$ is
$$
\hat G = \left\{ (a_0, \dots, a_7) \in \Z_3^8 \mid \sum_i a_i \equiv 0 \pmod 3
         \right\},
$$
where, if $g = (\zeta_0, \dots, \zeta_7) \in G$ and
$\alpha = (a_0, \dots, a_7) \in \hat G$, then
$\alpha(g) = \prod_i \zeta_i^{a_i}$.
To describe the action of $G$ on $H^6(\Hy, \C)$, we
use Griffiths' residue theory to represent cohomology classes with
differential forms with poles along $\Hy$ \cite{8}.

 Let
$\Omega = \sum_i (-1)^i x_i dx_0 \wedge \dots \wedge \widehat{dx_i} \wedge
\dots
\wedge dx_7$, where\ \ $\widehat{}$\ \ means omitted, be the rational top form
on~$\P^7$. Then the isomorphism $R^{3q-8}(\Hy) \to
H^{6-q+1,q-1}(\Hy)^{\text{prim}}$  is given by
$P \mapsto \Res\dfrac{P\Omega}{F^q}$. Since the Jacobian ideal of~$F$ is
generated by the monomials $x_0^2, \dots, x_7^2$,  the monomials of
degree~$i$ which do not contain squares form a basis for $R^i$;  moreover
$F$ is invariant under~$G$ and the
action of~$G$ commutes with~$\Res$, and hence
these monomials give a basis of eigenvectors. More precisely, the monomial $P =
\prod_i x_i^{b_i}$ gives an  eigenvector of eigenvalue $(b_0+1, \dots, b_7+1)$
(we are identifying the monomial~$P$ with the cohomogy class it corresponds.)
{}From what we have said, $b_i$ can only be~0 or~1, and hence the
characters with non--vanishing eigenspaces  are $\{ (a_0, \dots, a_7)
\in \hat G \mid    a_i \ne 0\ \forall i \}$. It is also clear from this
description that the eigenspaces are all one--dimensional. Moreover, since
$\Res\dfrac{P\Omega}{F^{q+1}}\in H^{6-q,q}(\Hy)$,
  we  see that the characters occurring in
$H^{6-q,q}(\Hy)$ are those for which $\sum_i a_i = 3(q+1)$, thinking of
the~$a_i$
as integers between~1 and~2.

In fact, these  eigenspaces are defined over $\Q(\zeta)$, and the action of
$\Gal(\Q(\zeta)/\Q) = \Z_3^*$ is given by $$
\sigma_a \cdot (a_0, \dots, a_7) = (aa_0, \dots, aa_7),  \tag3
$$
where $\sigma_a(\zeta) = \zeta^a$. This description allows us to analize the
$\Q$--structure of $H^6(\Hy, \C)$.

The cohomology class given by a monomial is only in $H^6(\Hy, \C)$, but since
the eigenspaces are one--dimensional, an appropriate multiple of an
eigenvector (by a transcendental number, usually) will be in $H^6(\Hy,
\Q(\zeta))$. These numbers are of course related to the periods of the Fermat
hypersurface, but since we don't need them explicitely, we will ignore their
precise values in what follows (also, we don't know these values; the periods
for Fermat hypersurfaces have been computed only for the cases of curves by
Rohrlich \cite{10, Appendix} and surfaces by Tretkoff \cite{20}.) We will use
the
following notation: if $\alpha\in\hat G$ is a character whose eigenspace is
nonzero, we let $P_{\alpha}$ be a monomial such that
$\Res\dfrac{P_{\alpha}\Omega}{F^4} \in H^6(\Hy, \Q(\zeta))$ is an eigenvector
relative to~$\alpha$.

We look now for rational cohomology classes in $H^{3,3}(\Hy)^{\text{prim}}$. By
the above discussion, a cohomology class in $H^6(\Hy, \Q(\zeta))$ is rational
if
and only if it is invariant under the action of the Galois group.
Let $\alpha\in\hat G$ be such that $\sum_i a_i = 12$, i.e., four of the $a_i$
are equal to~1 and the other four are equal to~2. Then $P_{\alpha}$ belongs to
$H^{3,3}(\Hy)^{\text{prim}}$ and the class corresponding to $P_{\alpha} +
P_{2\alpha}$ will be a primitive rational class. In this way we can construct
many rational primitive classes.

We now make explicit choices for the polynomials appearing in the formula of
Proposition~3.2, and compute the infinitesimal invariant of the corresponding
normal functions.
 Let $\alpha = (2,2,2,2,1,1,1,1)$ and $\beta =
(2,2,2,1,2,1,1,1)$, and let $a, b\in \Z$. Then $P_{\alpha}= Ax_0x_1x_2x_3$,
$P_{\beta} = Bx_0x_1x_2x_4$,
$P_{2\alpha} = Cx_4x_5x_6x_7$, $P_{2\beta} = Dx_3x_5x_6x_7$, where~$A$,
$B$, $C$ and~$D$ are nonzero complex numbers.
Let $\cy_{a,b}$ be
the primitive rational class corresponding to $P_{a,b} = a(P_{\alpha} +
P_{2\alpha}) + b(P_{\beta} + P_{2\beta})$.
Let $e = x_0x_1 + \dfrac{1}{C}x_2x_3 + \dfrac{1}{A}x_4x_5 +
x_6x_7 + \dfrac{h}{D}x_3x_5$, where $h$ is a transcendental number.
Let $Q = \dfrac{1}{A}x_4x_5x_6$ and
$R = \dfrac{1}{B}x_3x_5x_7$ be in $R^3(\Hy)$, so that
$Q\otimes R\in \ker\mu_{\Hy}$.

\proclaim{Proposition 4.1} Let the polynomials $P_{\alpha}$, $P_{2\alpha}$,
$P_{\beta}$,
$P_{2\beta}$, $e$, $Q$ and $R$ be chosen as above, and let $a$, $b\in\Z$.
Then:
\roster
\item if $a\ne0$, $a^2AC - b^2BD \ne0$,
the multiplication by $(P_{a,b}\cdot e)$ induces
      an isomorphism $f :R^1(\Hy) \to R^7(\Hy)$;
\item $\delta\nu_{\cy_{a,b}}(Q\otimes R) =  P_{a,b} Q (f^{-1}
(P_{a,b}R)) =  \dfrac{ab}{a+bh}\cdot x_0x_1x_2x_3x_4x_5x_6x_7 $.
\endroster
\endproclaim

\demo{Proof}
Let $\{f_0, \dots, f_7\}$ be a basis of~$R^1(\Hy)$, where $f_i$ is the image
of~$x_i$, and let $\{g_0, \dots, g_7\}$ be a basis of~$R^7(\Hy)$, where $g_i$
is the image of $\prod_{j\ne i}x_j$. It is then straightforward to check that,
using these bases, the matrix of the multiplication by $(P_{a,b}\cdot e)$ is
$$
M =
\pmatrix
0  &  a  &  0   &  0   &  0   &  0   &  0   &  0   \\
a  &  0  &  0   &  0   &  0   &  0   &  0   &  0   \\
0  &  0  &  0   &  aC  &  bD  &  0   &  0   &  0   \\
0  &  0  &  aC  &  0   &  0   &  bB  &  0   &  0   \\
0  &  0  &  bD  &  0   &  0   &  aA  &  0   &  0   \\
0  &  0  &  0   &  bB  &  aA  &  0   &  0   &  0   \\
0  &  0  &  0   &  0   &  0   &  0   &  0   & a+bh \\
0  &  0  &  0   &  0   &  0   &  0   & a+bh &  0
\endpmatrix
$$
Since $\det M = a^2 (a+bh)^2 (a^2AC-b^2BD)^2$, the multiplication by
$(P_{a,b}\cdot e)$ is an isomorphism as claimed.

A simple computation now gives that $P_{a,b}R = bg_6$, and
$f^{-1}(bg_6) = \dfrac{b}{a+bh}f_7$. The formula~(2) follows.
\qed \enddemo

We can now prove the theorem announced in the introduction:

\proclaim{Theorem 4.2} For the generic smooth cubic sevenfold, the Griffiths
group is not finitely generated.
\endproclaim

\demo{Proof}  The argument in \cite{21, 3.7} goes through word by
word, provided we show that the numbers $\dfrac{ab}{a+bh}$ span an infinite
dimensional $\Q$--vector space. Let then
$\sum_{i=1}^n c_i \dfrac{a_ib_i}{a_i+b_ih} = 0$ be a relation of linear
dependence over~$\Q$. We assume that $b_i=1$ for all~$i$, and that the $a_i$
are all distinct and non--zero. Clearing denominators, we obtain an
algebraic equation of degree~$n-1$ of which $h$ is a root. Since $h$ is
transcendental, all coefficients must be zero. This gives a homogeneous linear
system in the $c_i$'s, whose coefficient matrix has determinant a Van der Monde
in the $a_i$'s, and hence is non--zero. This implies that $c_i=0$ for all~$i$.
\qed
\enddemo


\Refs


\ref
\key 1 \by A. Albano
\paper Infinite generation of the Griffiths group: a local proof
\paperinfo Thesis, University of Utah, 1986
\endref

\ref
\key 2 \by F. Bardelli
\paper Curves of genus three on a general abelian threefold and the non--finite
generation of the Griffiths group
\inbook Arithmetic of complex manifolds (Erlangen, 1988)
\eds W. -P. Barth and H. Lange
\publ Springer--Verlag Lecture Notes in Mathemathics 1399
\publaddr New York Berlin Heidelberg
\yr 1989   \pages 10--26
\endref

\ref
\key 3 \by P. Candelas, E. Derrick, L. Parkes
\paper Generalized Calabi--Yau Manifolds and the mirror of a rigid manifold
\finalinfo preprint CERN-TH.6831/93, UTTG-24-92
\endref


\ref
\key 4 \by C. Ciliberto, J. Harris, R. Miranda
\paper General components of the Noether--Lefschetz locus and their density
       in the space of all surfaces
\jour  Math. Ann.  \vol 282  \yr 1988  \pages 667--680
\endref

\ref
\key 5 \by H. Clemens
\paper Homological equivalence, modulo algebraic equivalence, is not
finitely generated
\jour Inst. Hautes \'Etudes Sci. Publ. Math.
\vol 58  \yr 1983  \pages 19--38
\endref


\ref
\key 6 \by P. Deligne
\paper Cohomologie des intersection compl\`etes
\inbook Groupe de Monodromie en G\'eom\'etrie Alg\'ebrique (\rom{S.G.A. 7 II})
\publ Springer--Verlag Lecture Notes in Mathemathics 340
\publaddr New York Berlin Heidelberg
\yr 1973   \pages 39--61
\endref


\ref
\key 7 \by M. Green
\paper Griffiths' infinitesimal invariant and the Abel--Jacobi map
\jour Jour. Diff. Geom.  \vol 29  \yr 1989  \pages 545--555
\endref



\ref
\key 8 \by P. A. Griffiths
\paper On the periods of certain rational integrals. I, II
\jour Ann. of Math.  \vol 90  \yr 1969  \pages 460--541
\endref


\ref
\key 9 \bysame
\paper Infinitesimal invariant of normal functions
\inbook Topics in transcendental algebraic geometry
\ed P. A. Griffiths
\bookinfo Annals of Mathematics Studies 106
\publ Princeton University Press \publaddr Princeton
\yr 1984\pages 305--316
\endref




\ref
\key 10 \by B. H. Gross
\paper  On the periods of abelian integrals and a formula of Chowla and
        Selberg (with an appendix by D. E. Rohrlich)
\jour Inventiones Math.  \vol 45  \yr 1978  \pages 193--211
\endref


\ref
\key 11 \by S. Kleiman
\paper Geometry on Grassmannians and applications to splitting bundles and
       smoothing cycles
\jour Publ. Math. I.H.E.S.  \vol 36  \yr 1969  \pages 281--298
\endref

\ref
\key 12 \by D. Morrison
\paper Mirror symmetry and rational curves on quintic threefolds: A guide
for the mathematicians
\jour Jour. A.M.S.  \vol 6  \yr 1993  \pages 223--247
\endref


\ref
\key 13 \by M. Nori
\paper Cycles on the generic abelian threefold
\jour Proc. Indian Acad. Sci.  \vol 99  \yr 1989  \pages 191--196
\endref

\ref
\key 14 \bysame
\paper Algebraic cycles and Hodge theoretic connectivity
\jour Inventiones Math. \vol 111  \yr 1993  \pages 349--373
\endref


\ref
\key 15 \by A. Ogus
\paper Griffiths transversality and crystalline cohomology
\jour Ann. of Math.  \vol 108  \yr 1978  \pages 395--419
\endref

\ref
\key 16 \by Z. Ran
\paper Cycles on Fermat hypersurfaces
\jour Compositio Math.  \vol 42  \yr 1980/81  \pages 121--142
\endref


\ref
\key 17 \by C. Schoen
\paper Complex multiplication cycles on elliptic modular threefolds
\jour Duke Math. Jour.  \vol 53  \yr 1986  \pages 771--794
\endref




\ref
\key 18 \by T. Shioda
\paper The Hodge conjecture for Fermat varieties
\jour Math. Ann. \vol 245  \yr 1979  \pages 175--184
\endref


\ref
\key 19 \by J. Steenbrink
\paper Some remarks about the Hodge conjecture
\inbook Hodge theory  \eds E. Cattani, F. Guill\'en, A. Kaplan and F. Puerta
\publ Springer--Verlag Lecture Notes in Mathemathics 1246
\publaddr New York Berlin Heidelberg
\yr 1987   \pages 165--175
\endref

\ref
\key 20 \by M. D. Tretkoff
\paper The Fermat surface and its periods
\inbook Recent developments in several complex variables
\ed J. E. Fornaess
\bookinfo Annals of Mathematics Studies 100
\publ Princeton University Press  \publaddr Princeton
\pages 413--428 \yr 1981
\endref


\ref
\key 21 \by C. Voisin
\paper Une approche infinit\'esimale du th\'eor\`eme de H. Clemens sur les
       cycles d'une quintique g\'en\'erale de $\P^4$
\jour Jour. of Alg. Geom. \vol 1  \yr 1992  \pages 157--174
\endref


\ref
\key 22 \bysame
\paper Sur l'application d'Abel--Jacobi des vari\'et\'es de Calabi--Yau
de dimension trois
\jour Ann. Sci. E.N.S.  \toappear
\endref






\endRefs
\enddocument